\documentclass{amsart}

\usepackage[english]{babel}

\usepackage[normalem]{ulem} 
\usepackage[super]{nth} 
\usepackage{setspace}
\usepackage{amsmath, amsthm,amssymb}
\usepackage{float}
\usepackage{bm}
\usepackage{graphicx}
\usepackage{pgfplots}
\usepackage{pgfplotstable}
\pgfplotsset{compat = newest}

\renewcommand{\le}{\leqslant}

\renewcommand{\leq}{\leqslant}

\begin{document}

\title{An Exact Closed-Form Solution of the Lotka-Volterra Equations}

\author{Jean-Luc Boulnois}

\address{Babson College, Babson Park, Wellesley, Massachusetts 02457, jlboulnois@msn.com}

\begin{abstract}
The classical Lotka-Volterra predator-prey system is often used in species competition modeling. An exact, closed-form solution is derived when the natural growth rate of the prey species and decay rate of the predators are equal in magnitude. A standard functional  transformation yields a novel system of two \textit{partially uncoupled} first-order ODEs for ``hybrid-species'', with one being \textit{autonomous}. New exact, closed-form time-dependent solutions are derived for each individual species. An analytical expression for the system's oscillation period valid for any value of the system's energy is  derived in terms of a novel universal function. 
\end{abstract}


\maketitle


\section{Introduction}

The historic Lotka-Volterra (``LV'') predator-prey  system of two coupled first-order nonlinear differential equations has first been investigated in ecological and chemical systems \cite{volterra1926}, \cite{lotka}. This idealized model describes the competition of two isolated coexisting species: a `prey' population evolves while feeding from an infinitely large resource supply, whereas `predators' interact by exclusively feeding on preys, either through direct predation or as parasites. As a result the respective populations exhibit undamped oscillations as a function of time with a period which depends on the species interaction rates.
\newline

The  classical LV model is based on four time-independent, positive, and constant rates with two  representing species self-interaction, i.e. natural exponential growth rate $\alpha$ and decay rate $\delta$ per individual of the respective prey and predator populations, and two others  characterizing inter-species interaction.

\section{Normalized Equations }

Without any loss of generality, the LV system of two coupled first order ordinary differential equations (ODE) can  be simplified by simultaneously scaling the predator and prey populations together with time through a dimensionless time $t$ based on the factor $1/\sqrt{\alpha \delta}$. The system is shown to only depend on a positive coupling parameter $\lambda $, ratio of the respective growth and decay rates of each species taken separately, defined as
\begin{equation} \label{EQ_1_} 
\lambda =\sqrt{\frac{\alpha }{\delta} }   
\end{equation} 

The respective instantaneous populations of preys and predators, labeled $u$ and $v$, both $\geqslant 0$, are assumed to be continuous functions of time: a normalized form of the LV system is obtained as a set of two coupled first-order nonlinear ODEs solely depending on this single coupling ratio  $\lambda $, \cite{boulnois}.
\newline

We consider the special case when the prey population natural growth rate and predator population natural decay rate have \textit{equal magnitude and sign}, i.e.  when $\lambda = 1$. An exact solution has previously been derived \cite{varma} in the  case when the prey growth rate and predator decay rate are identical in magnitude, but with \textit{opposite signs}, i.e. $\alpha = - \delta$, a condition precluding population oscillation. 
\newline

In this special case, the two species time-evolution is  modeled as a system of two coupled  autonomous nonlinear ODEs where the ``dot'' on  $\dot{u}$  and  $\dot{v}$  indicates a derivative with respect to the  time  $t$
\begin{subequations}\label{EQ_2_}
    \begin{align}
        \dot{u}&= u\left(1-v\right) \quad     \text{for preys} \\
        \dot{v}&= v\left(u-1\right) \quad     \text{ for predators}
    \end{align}
\end{subequations}

Numerous solutions of system \eqref{EQ_2_} have been developed including trigonometric series \cite{frame}, mathematical transformations \cite{evans1999}, Taylor series expansions \cite{Scarpello}, perturbation techniques \cite{rao}, and Lambert W-functions \cite{shih2005}. 
\newline

The system \eqref{EQ_2_} is non-trivial but is known to possess a dynamical invariant or ``constant of motion'' representing the conservation of the positive, constant energy ``$h$'' of the system, which can be written

\begin{equation} \label{EQ_3_} 
h = \frac{1}{2} \left (u+v-\ln (u v) \right)-1
\end{equation} 

In the following sections, through an exponential functional  transformation  we introduce ``hybrid species" within a new set of two \textit{partially uncoupled} first-order ODEs with one being \textit{autonomous}. A new, exact closed-form solution is  derived for each hybrid-species  separately in terms of a single  quadrature. An exact analytical solution of the LV system for each individual prey and predator species $u(t)$ and $v(t)$ is derived as a function of time. The exact population oscillation period is further presented in terms of a novel universal energy function.

\section{Exact Solutions with Hybrid Predator-Prey Species}

Upon introducing a functional transformation of the original prey and predator population set $\{$\textit{u(t), v(t)}$\}$ to a new set of ``\textit{hybrid} predator-prey species" $\{$\textit{$\xi$(t)}, \textit{$\eta$(t)}$\}$  with $\xi$  and $\eta \in \mathbb{R}$, representing the symbiotic coupling between interacting preys and predators according to 
\begin{subequations} \label{EQ_4_} 
\begin{align}
    u(t)  &= e^{\xi(t)-\eta(t)}     \quad \text{for preys }  \label{EQ_4_a}\\
    v(t)  &= e^{\xi(t)+\eta(t)}   \quad \text{for predators } \label{EQ_4_b}
     \end{align} 
\end{subequations} 

the energy conservation equation \eqref{EQ_3_} becomes
\begin{equation} \label{EQ_5_} 
h = e^{\xi } \cosh (\eta ) -\xi -1 
\end{equation} 

This energy conservation relationship \eqref{EQ_5_}  is recast into a  form which provides a natural separation  between the  functions  \textit{$\eta$(t)} and \textit{ $\xi$(t)}, with the system's positive energy $h$ explicitly associated with the $\xi$-function only \cite{boulnois}
\begin{equation} \label{EQ_6_} 
\cosh (\eta ) =(h+1+\xi )e^{-\xi }  
\end{equation} 

Equation \eqref{EQ_6_} can equivalently be written in terms of the inverse hyperbolic cosine function
\begin{equation} \label{EQ_7_} 
\eta ^{\pm } (\xi )=\pm \cosh ^{-1} \bigl((h+1+\xi )e^{-\xi } \bigr)
\end{equation} 

In the following we define a useful compact  auxiliary function  \textit{U($\xi$)} $\geqslant 1$ appearing throughout  
\begin{equation} \label{EQ_8_} 
U(\xi )=(h+1+\xi )e^{-\xi }  
\end{equation} 

The hybrid-species population $\xi (t)$  thus oscillates between the respective negative and positive roots $\xi ^{-} (h)$ and  $\xi ^{+} (h)$, solutions of the equation $U(\xi )=1$ displayed in Table \ref{table:1} below for several increasing values of $h$
\begin{equation} \label{EQ_9_} 
e^{\xi } -\xi -1=h  \quad \text{ with } h\geqslant  0
\end{equation} 

For any value of the constant energy $h$, in the $\xi -\eta $ phase plane,  Eq. \eqref{EQ_6_} represents a  closed-orbit loop consisting of two respective symmetric branches  $\eta ^{+} (\xi )$ and  $\eta ^{-} (\xi )$ around the fixed point $(0,0)$. This orbit is bounded by the limits $\xi ^{-} (h)$ and  $\xi ^{+} (h)$ on the $\eta =0$ horizontal axis, and since  \textit{U($\xi$)}  admits a maximum $e^{h} $ located at $\xi = -h$ , it is also bounded vertically by the two respective roots  of the equation  $\eta ^{\pm } (-h)= \pm \cosh^{-1} (e^{h} )$. 
\newline

Lastly, upon inserting the exponential transformation \eqref{EQ_4_} into the normalized LV system \eqref{EQ_2_}, a new partially uncoupled system of two \nth{1} order ODEs is obtained for each hybrid species taken separately,  with  the  evolution  of the $\xi$ hybrid species  population  represented by a  \nth{1} order nonlinear \textit{autonomous} ODE
\begin{subequations}\label{EQ_10_}
    \begin{align}
        \dot{\eta }&=\xi + h \label{EQ_10_a}\\
       \dot{\xi }&= \pm \bigl( (h+1+\xi )^{2} -e^{2\xi } \bigr)
        ^{1/2} \label{EQ_10_b}
    \end{align}
\end{subequations}

The solution of   system \eqref{EQ_10_} represents the  time-evolution of the hybrid-species $\eta (t)$  and $\xi (t)$, \textit{albeit} due to the hybrid species transformation \eqref{EQ_4_}, equation \eqref{EQ_10_a} becomes \textit{linear}. Remarkably, in this $\lambda = 1$ case, the solution of the LV problem is considerably simplified since  a solution of  the autonomous equation \eqref{EQ_10_b}  only is required.
\newline

A numerical solution for $\xi (t)$ can be obtained by integrating Eq. \eqref{EQ_10_b} using a standard fourth-order Runge-Kutta (RK4) method. Figure 1, reprinted from \cite{boulnois}, presents the $\xi (t)$-solution obtained by numerical integration of  \eqref{EQ_10_b} for an energy $h = 2$ with initial condition $\xi(0) = \xi ^{-} (h)$. The growth and decay phases of the even function $\xi (t)$ are observed to be symmetric  relative to the half-period  $t^{*}$ when $\xi (t^{*}) = \xi^{+}(h)$. The LV system  solution is finalized for the two branches  $\eta ^{\pm } (t)$  by inserting $\xi (t)$ derived above into Eq. \eqref{EQ_7_} since $\eta(0) = 0$.
\newline

\begin{figure}[H]
  \centering
 \includegraphics{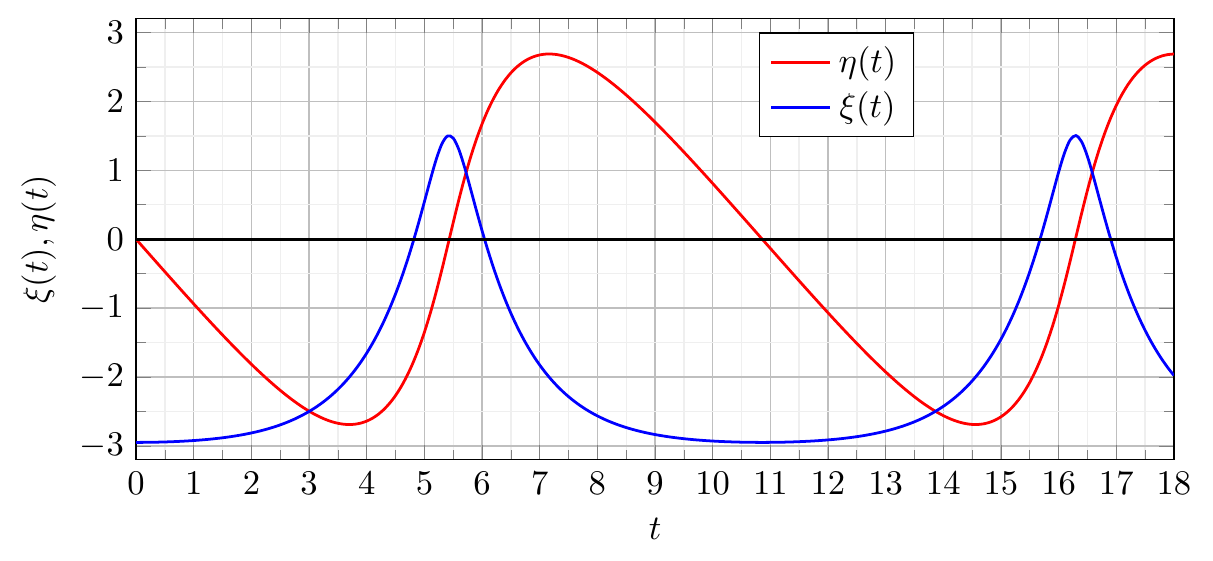}

  \caption{Solutions for $\xi (t)$ and  $\eta (t)$ as a function of time $t$  obtained by numerical integration of the quadrature solution Eq. \eqref{EQ_10_b} with energy $h = 2$ (from \cite{boulnois})}.
  \label{fig:Fig1}

\end{figure}

Over the respective intervals $\xi ^{-} \le \xi(t) \le \xi ^{+} $ and $\xi ^{+} \geqslant  \xi(t) \geqslant  \xi ^{-} $  corresponding to the growth and decay phases of the hybrid species population $\xi(t)$, an integral expression for  $t(\xi )$  is readily obtained by performing the integration with the respective positive root (growth phase) and negative root (decay phase) in \eqref{EQ_10_b}, yielding the following  \textit{quadrature} solution 
\begin{subequations} \label{EQ_11_} 
   \begin{align}
        t(\xi )=\int _{\xi ^{-} }^{\xi }\frac{dx}{\sqrt{(h+1+x)^{2} - e^{2x} } }  \label{EQ_11_a}\\ 
        t(\xi )= t^{*} + \int _{\xi ^{+} }^{\xi }\frac{dx}{\sqrt{(h+1+x)^{2} - e^{2x} } }  \label{EQ_11_b}
   \end{align}
\end{subequations} 

Even though the $\xi (t)$ hybrid-species population oscillation is not explicitly expressed as a function of time $t$, the function  \textit{t($\xi$)}  being monotonic and continuous on each respective integration interval,  its inverse function $t^{-1}: \mathbb{R} \to \mathbb{R}$ defined by $\xi (t) = t^{-1}(\xi)$, which only  depends on the energy level $h$, exists and is unique, monotonic, and continuous on each interval. At the respective limits $\xi ^{-} (h)$ and $\xi ^{+} (h)$  the integrand of  \eqref{EQ_11_} has a weak singularity of the square root type, but is strictly continuous  over the interval and the integral is convergent.
\newline

Together with \eqref{EQ_7_}, the exact solution \eqref{EQ_11_} for $\xi(t)$ over the respective intervals $\xi ^{-} \le \xi \le \xi ^{+} $ and $\xi ^{+} \geqslant  \xi \geqslant  \xi ^{-} $ constitutes the final solution of the LV problem for the ``hybrid species" in the special $\lambda = 1$ case considered here.
\newline

The  solution \eqref{EQ_11_} is similar in form to a solution derived by Evans and Findley (Eq. (17) in \cite{evans1999}); however, the above integral expression lends itself to simpler analytical or numerical integration. An exact expression for \eqref{EQ_11_} is further proposed in Appendix 1 in terms of exponential integral functions. 

\section{Exact Solutions for the Prey and Predator Species Populations}

Exact solutions for the time evolution of the prey and predator populations are derived by inserting the respective hybrid-species populations $\xi (t)$  and  $\eta (t)$  obtained from Eqs. \eqref{EQ_11_} and  \eqref{EQ_7_} into the original definition \eqref{EQ_4_}.  This results  in two \textit{uncoupled} solutions for the individual populations $u(t)$ and $v(t)$ of the prey and predator species. Over the growth and decay phases of the symmetric $\xi(t)$ function, these exact \textit{uncoupled} analytical solutions are expressed as
\newline

\underline{Interval $0 \le t \le t^*$} $\mapsto$ interval $\xi ^{-} \le \xi(t) \le \xi ^{+} $, i.e. $\xi (t)$ growth phase 
\begin{subequations}\label{EQ_12_}
    \begin{align}
      u(t)&= h+1+\xi (t) +  \sqrt{(h+1+\xi(t))^{2} - e^{2\xi(t)} }   \quad &\text{for preys } \label{EQ_12_a}\\
      v(t)&=  h+1+\xi (t) -  \sqrt{(h+1+\xi(t))^{2} - e^{2\xi(t)} }   \quad &\text{for predators }\label{EQ_12_b}
    \end{align}
\end{subequations}

 with $\xi (t)= t^{-1}(\xi)$  derived from \eqref{EQ_11_a} together with $\xi(0) = \xi ^{-}(h)$.
\newline

\underline{Interval $t^{*} \le t \le 2t^{*}$} $\mapsto$ interval $\xi ^{+} \geqslant  \xi(t) \geqslant  \xi ^{-} $, i.e. $\xi (t)$ decay phase 
\begin{subequations}\label{EQ_13_}
    \begin{align}
      u(t)&=  h+1+\xi (t) -  \sqrt{(h+1+\xi(t))^{2} - e^{2\xi(t)} }   \quad &\text{for preys } \label{EQ_13_a}\\
      v(t)&=  h+1+\xi (t) +  \sqrt{(h+1+\xi(t))^{2} - e^{2\xi(t)} }   \quad &\text{for predators }\label{EQ_13_b}
    \end{align}
\end{subequations}

with $\xi (t)=t^{-1}(\xi)$  derived from  \eqref{EQ_11_b} together with $\xi(t^{*}) = \xi ^{+}(h)$.
\newline

Figure 2 displays the  exact uncoupled analytical solutions for the time evolution of the preys $u(t)$ and the predators $v(t)$ when their respective growth and decay rates are equal in magnitude, and when the system's energy is $h = 2$.
\newline

It is observed that the prey population $u(t)$ exhibits an initial growth at rate significantly slower than its own fast decay rate while the opposite is the case for the predators; also the peak population of the preys occurs when the population is mature, i.e. when the predator population is young $v=1$, and vice-versa. 

\begin{figure}
  \centering
  \includegraphics{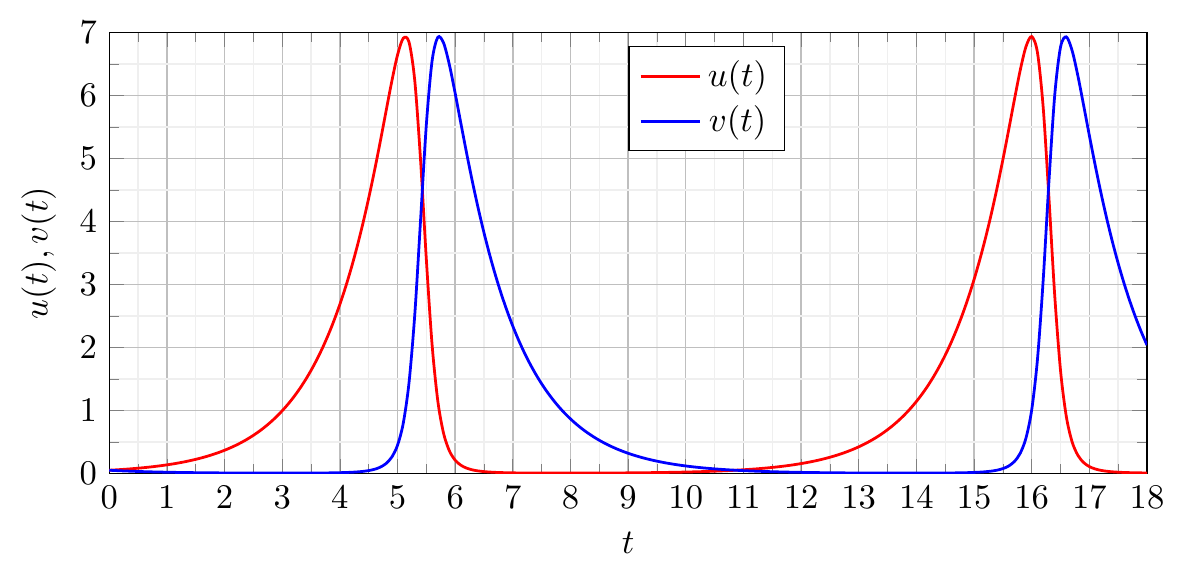}
  
  \caption{Exact analytical solutions for $u(t)$ and  $v(t)$ as a function of time $t$  obtained from Eqs. \eqref{EQ_12_} and \eqref{EQ_13_} with energy $h = 2$ }.
  \label{fig:Fig2}
\end{figure}

\section{Oscillation Period of the LV System}

In the special  case  $\lambda =1$ considered here, when the rates $\alpha$ and $\delta$ are equal, the exact LV system period $T(h)$, which has been shown to be the shortest for any energy $``h"$ \cite{boulnois}, can uniquely be expressed in terms of a universal energy function  $\Theta (h)$  as
\begin{equation} \label{EQ_14_} 
T(h)=\frac{2\pi }{\alpha  } \Theta (h) 
\end{equation} 

The  energy function $\Theta (h)$ introduced in \cite{boulnois} is defined by integrating \eqref{EQ_11_} over the entire $\xi$-interval 
\begin{equation} \label{EQ_15_} 
\Theta (h)=\frac{1}{\pi } \int _{\xi ^{-} }^{\xi ^{+} }\frac{dx}{\sqrt{(h+1+x)^{2} - e^{2x}} }   
\end{equation} 

At small orbital energy ($h \ll 1$) where $\xi^{\pm} (h) = \pm \sqrt{2h}$,  the function $\Theta (h)$ is directly expressed in terms of the complete elliptic integral of the first kind $\bm{K}(k)$ with  modulus $k$ 
\begin{equation} \label{EQ_16_}
\Theta (h) = \frac{1}{\sqrt{1+\sqrt{2h}}}\frac{2}{\pi} \bm{K}(k) \quad \text{ with } \quad   k = \sqrt{\frac{2\sqrt{2h}}{1+\sqrt{2h}}}
\end{equation}

A standard series expansion for $\bm{K}(k)$ yields
\begin{equation} \label{EQ_17_}
\Theta (h) = 1+\frac{1}{3} h + \frac {1}{42} h^2 + O(h^3)
\end{equation}

For small oscillation amplitudes,  the integral \eqref{EQ_15_} becomes  independent of the  energy $h$  and is exactly equal to $\pi$, hence $\Theta(h) = 1$ in \eqref{EQ_14_}; the LV system becomes that of two  \textit{coupled harmonic oscillators} for which the period $T(h)$  solely depends on the pulsation $\alpha $, as already established  \cite{volterra1926}, \cite{waldvogel}.
\newline

At high orbital energy ($h \gg 1$), the contribution from the exponential term in \eqref{EQ_15_} becomes negligible since $\xi < 0$ over most of the integration interval except when $\xi$ approaches $\xi^{+} (h)$: since by definition $\xi(t) \geqslant   \xi ^{-} (h)$, approximating the exponential term by its lowest value  $e^{2\xi ^{-} (h)} $ and performing the integration yields a useful asymptotic expression for $\Theta(h)$
\begin{equation} \label{EQ_18_} 
\Theta _{\text{asymp}} (h)\cong \frac{1}{\pi } \left( \xi ^{+} (h)-\xi ^{-} (h) + ln(2) \right)    \quad   \text{  with } h \gg 1
\end{equation} 

\begin{table}[H]
{\small

\begin{tabular}{|p{0.38in}|p{0.38in}|p{0.38in}|p{0.38in}|p{0.38in}|p{0.38in}|p{0.38in}|p{0.38in}|p{0.38in}|} \hline 
\textit{h} & 0.3 & 0.5 & 1 & 2 & 3 & 5 & 7 & 10 \\ \hline 
$\xi ^{-} (h)$ & -0.889 & -1.198 & -1.841 & -2.948 & -3.981 & -5.998 & -8.000 & -11.00 \\ \hline 
$\xi ^{+} (h)$ & 0.686 & 0.858 & 1.146 & 1.505 & 1.749 & 2.091 & 2.336 & 2.611 \\ \hline
$\Theta (h)$ & 1.102 & 1.173 & 1.355 & 1.728 & 2.102 & 2.828 & 3.535 & 4.569  \\ \hline
\end{tabular}
}
\caption{Roots of  $e^{\xi } -\xi -1=h$  as a function of the energy $h$   from Eq. \eqref{EQ_9_}, and values of $\Theta (h)$ from \eqref{EQ_15_}}
\label{table:1}
\end{table}

The universal function $\Theta (h)$ tabulated in Table 1 for various values of the energy $h$  is also displayed on Fig. 3 : $\Theta (h)$, and by extension $ T(h)$,  is a monotonically increasing function of the energy-dependent amplitude  $\xi ^{+} (h)-\xi ^{-} (h)$ of the $\xi(t)$ function only  \cite{waldvogel}. Also shown is the asymptotic approximation \eqref{EQ_18_} which is practically indistinguishable from the  exact function  $\Theta(h)$  for $h \geqslant  4$.

\begin{figure}[H]
  \centering
  \includegraphics{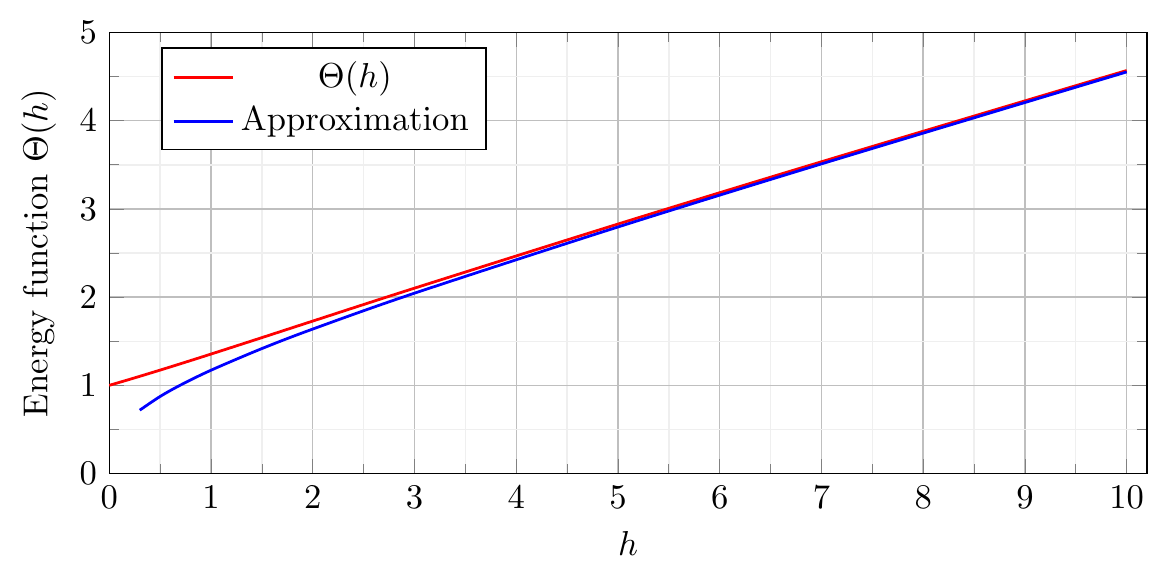}
  
  \caption{Universal energy function $\Theta(h)$ \eqref{EQ_15_} and asymptotic approximation \eqref{EQ_18_} for $h \gg 1$}
  \label{fig:Fig3}
\end{figure}

Upon comparing the methods of Volterra \cite{volterra1926}, Hsu \cite{hsu}, Waldvogel \cite{waldvogel}, and Rothe \cite{rothe}, Shih demonstrated that all of their integral representations for the period of the two-species LV system are equivalent to his own solution in terms of a sum of several convolution integrals \cite{shih1997}. 
The period derived here when $\lambda = 1$ is expressed as a single integral \eqref{EQ_15_}.

\section{Conclusion}

The coupled \nth{1} order non-linear ODE system for the LV problem of two interacting species has been analyzed in the special case when the relative growth/decay rates of each species taken independently are equal. 
Based on a standard functional  transformation introducing ``hybrid-species populations", a new  set of two \nth{1} order ODEs is obtained with one being \textit{autonomous}. 
\newline

In this special case, the LV problem partially uncouples and an exact explicit closed-form solution is derived in terms of the system's orbital energy  $h$ as a simple quadrature for the time evolution  of the hybrid-species population $\xi (t)$;  the other hybrid species' solution $\eta (t)$  is explicitly expressed in terms of the former (Eqs.\eqref{EQ_11_} and \eqref{EQ_7_}).
As a result, exact \textit{uncoupled} analytical solutions for each of the original prey and predator populations $u(t)$ and $v(t)$ are derived as a function of time.
\newline

Further, an exact, closed-form expression for the non-linear LV system oscillation period $T(h)$ is derived in terms of a universal LV energy function together with a simple asymptotic expression for high  energies $(h \gg 1)$.

\section*{Appendix 1}

 Upon recalling the definition \eqref{EQ_8_} of $U(\xi)$,  a series expansion for the quadrature solution \eqref{EQ_11_a} is derived by first expressing the integral as
\begin{equation*} \label{EQ_A2.1_} 
t(\xi )=\cosh ^{-1} \bigl(U(\xi )\bigr)+ e^{\xi} \sqrt {U(\xi)^{2}-1}  +\int _{\xi ^{-} }^{\xi } \frac{e^{2\xi}}{(h+1+\xi)\sqrt{1-U(x)^{-2} } }  dx               \tag{A2.1}
\end{equation*} 

Upon observing in Table 1 that as $h$ increases, $\xi^{-}(h) \simeq -(h+1)$, the exponential factor in the integral  thus becomes negligible for most of the integration interval  up to $\xi \simeq 0$.  Consequently the contribution to the solution  $t(\xi)$ in the growth phase $(\xi \leq 0)$ of the  $\xi$-function principally comes from the first two terms in \eqref{EQ_A2.1_}, while near its maximum $\xi ^{+}(h)$ where both of these terms vanish,  the contribution over the interval $\xi \simeq 0$ to $\xi^{+}(h)$  mostly comes from the integral.
\newline

Since $1\le U(\xi )\le e^{h} $, a binomial expansion of the integrand with  binomial coefficients expressed in terms of the Gamma function $\Gamma (p)$  yields an exact solution in terms of a converging series
 
\begin{equation*}  \label{EQ_A2.2_} 
t(\xi )=\cosh ^{-1} \bigl(U(\xi ) \bigr) +e^{\xi} \sqrt {U(\xi)^{2}-1} + \sum _{p=0}^{\infty }\frac{\Gamma \left(\frac{1}{2} \right)}{\Gamma \left(\frac{1}{2} -p \right)\Gamma \left (p+1 \right)} \int _{\xi ^{-} }^{\xi } \frac{e^{2x}U(x)^{-2p}}{h+1+x} dx                \tag{A2.2}
\end{equation*} 

The first integral $(p=0)$ is directly expressed in terms of the exponential integral function Ei(x)
\begin{equation*}  \label{EQ_A2.3_} 
\int _{\xi ^{-} }^{\xi } \frac{e^{2x}}{h+1+x} dx =  e^{-2(h+1)}\bigl(\mathrm{Ei}({2(h+1+\xi)}) -\mathrm{Ei}({2e^{2\xi^{-}}})\bigr)             \tag{A2.3}
\end{equation*} 

When inserted into \eqref{EQ_A2.2_} this expression provides a zeroth order $(p=0)$ solution for $t(\xi)$, hence for $\xi(t) = t^{-1}(\xi)$ as discussed earlier.
\newline

When the integer $p$ is $1, 2, 3, \dots$,  each integral $I_{2p}(\xi)$ in the expansion \eqref{EQ_A2.2_} is of the form 
\begin{equation*} \label{EQ_A2.4_} 
I_{2p} (\xi )=\int _{\xi ^{-} }^{\xi }\frac{e^{2px} dx}{(h+1+x)^{2p+1} }                              \tag{A2.4}
\end{equation*} 

Successive integration by parts and substitution into \eqref{EQ_A2.2_} result in a  convergent series of exponential integral functions with positive argument of the form $ e^{-2p(h+1)}\mathrm{Ei}\bigl({2p(h+1+\xi )}\bigr)$. 
\newline


\nocite{}

\bibliographystyle{abbrv}


\begin{thebibliography}{10}


\bibitem{boulnois}
J.~L. Boulnois.
\newblock {P}redator-{P}rey linear coupling with hybrid species.
\newblock {\em arXiv}, 2301.00673, 2022.






\bibitem{evans1999}
C.~M. Evans and G.~L. Findley.
\newblock A new transformation of the {L}otka-{V}olterra problem.
\newblock {\em J. Math. Chem.}, 25(Added Volume):105--110, 1999.

\bibitem{frame}
J.~Frame.
\newblock Explicit solutions in two species volterra systems.
\newblock {\em Journal of Theoretical Biology}, 43(1):73 -- 81, 1974.



\bibitem{hsu}
S.~B. Hsu.
\newblock A remark on the period of the periodic solution in the
  {L}otka-{V}olterra system.
\newblock {\em J. Math. Anal. Appl.}, 95(2):428--436, 1983.



\bibitem{lotka}
A.~J. Lotka.
\newblock Undamped oscillations derived from the law of mass action.
\newblock {\em Journal of the American Chemical Society}, 42(8):1595--1599,
  1920.



\bibitem{plank}
M.~Plank.
\newblock Hamiltonian structures for the {$n$}-dimensional {L}otka-{V}olterra
  equations.
\newblock {\em J. Math. Phys.}, 36(7):3520--3534, 1995.

\bibitem{rao}
D.~V.~G. Rao and Y.~L.~P. Thorani.
\newblock A study of the solutions of the {L}otka-{V}olterra prey-predator
  system using perturbation technique.
\newblock {\em Int. Math. Forum}, 5(53-56):2667--2673, 2010.

\bibitem{rothe}
F.~Rothe.
\newblock The periods of the {V}olterra-{L}otka system.
\newblock {\em J. Reine Angew. Math.}, 355:129--138, 1985.

\bibitem{Scarpello}
G.~Mingari~Scarpello and D.~Ritelli.
\newblock A new method for the explicit integration of {L}otka-{V}olterra
  equations.
\newblock 11:1--17, 01 2003.

\bibitem{shih1997}
S.-D. Shih.
\newblock The period of a {L}otka-{V}olterra system.
\newblock {\em Taiwanese J. Math.}, 1(4):451--470, 12 1997.

\bibitem{shih2005}
S.-D. Shih.
\newblock Comments on “a new method for the explicit integration of
  {L}otka-{V}olterra equations”.
\newblock {\em Divulgaciones Matemáticas}, 13(2):99--106, 2005.



\bibitem{varma}
V.~S. Varma.
\newblock Exact solutions for a special prey-predator or competing species
  system.
\newblock {\em Bull. Math. Biology}, 39(5):619--622, 1977.

\bibitem{volterra1926}
V.~Volterra.
\newblock {Variation and fluctuations of the number of individuals of animal
  species living together}.
\newblock In R.~N. Chapman, editor, {\em Animal Ecology}, pages 31--113.
  McGraw-Hill, 1926.

\bibitem{waldvogel}
J.~Waldvogel.
\newblock The period in the {L}otka-{V}olterra system is monotonic.
\newblock {\em J. Math. Anal. Appl.}, 114(1):178--184, 1986.





\end{thebibliography}




\end{document}